\documentclass[10pt,reqno]{amsart}
\usepackage{amsfonts,amsthm,latexsym,amsmath,amssymb,amscd,amsmath,
epsf}

\theoremstyle{plain}
\newtheorem{theorem}           {Theorem}

\newtheorem{corollary}         {Corollary}
\newtheorem{lemma}             {Lemma}

\theoremstyle{definition}
\newtheorem{definition}        {Definition}
\newtheorem{problem}           {Problem}

\newtheorem{notation}          {Notation}

\newtheorem*{ack}              {Acknowledgements}

\theoremstyle{remark}
\newtheorem{remark}            {Remark}

\renewcommand{\Im}{\mathfrak{Im}}
\renewcommand{\Re}{\mathfrak{Re}}
\newcommand{\al}{\alpha}
\newcommand{\be}{\beta}
\newcommand{\la}{\lambda}

\newcommand {\bC} {\mathbb C}
\newcommand {\bR} {\mathbb R}

\newcommand {\eps} {\epsilon}

\newcommand{\bN}{\mathbb N}
\newcommand{\HH}{\mathcal{H}}

\newcommand{\Om}{\Omega}

\newcommand{\NN}{\mathbb N}

\newcommand{\RR}{\mathbb{R}}
\newcommand{\CC}{\mathbb{C}}

\newcommand{\HC}{\overline{\mathcal{H}}}

\input cyracc.def
\newfam\cyrfam
\font\tencyr=wncyr10
\font\sevencyr=wncyr7

\textfont\cyrfam=\tencyr \scriptfont\cyrfam=\sevencyr
    \scriptscriptfont\cyrfam=\sevencyr

\font\tencyr=wncyr10

\def\newop#1{\expandafter\def\csname #1\endcsname{\mathop{\rm
#1}\nolimits}}

\newop{slc}
\newop{sign}
\newop{mdeg}

\begin{document}
\numberwithin{equation}{section}
\title[P\'olya-Schur Master Theorems]{P\'olya-Schur Master Theorems for 
Circular Domains and their Boundaries}
\author[J.~Borcea]{Julius Borcea}
\address{Department of Mathematics, Stockholm University, SE-106 91
Stockholm, Sweden}
\email{julius@math.su.se}
\author[P.~Br\"and\'en]{Petter Br\"and\'en}
\address{Department of Mathematics, Royal Institute of Technology,
SE-100 44 Stockholm, Sweden}
\email{pbranden@math.kth.se}



%\thanks{$^*$ Corresponding author.}

\begin{abstract}
We characterize all linear operators on finite or infinite-dimensio\-nal
polynomial spaces that preserve the property of 
having the zero set inside a prescribed region $\Om\subseteq \bC$ for
arbitrary closed 
circular domains $\Om$ (i.e., images of the closed unit disk under a 
M\"obius transformation) and their boundaries. This provides
a natural framework for dealing with several long-standing fundamental 
problems, which we solve in a unified way. In particular, for $\Om=\bR$ our 
results settle open questions that go back to Laguerre and P\'olya-Schur. 
\end{abstract}

\maketitle

\section{Introduction}\label{s1}

Some of the main challenges in the theory of distribution of zeros of 
polynomials and transcendental entire functions concern the description of 
linear operators that preserve certain prescribed (``good'') properties.
Notwithstanding their fundamental character, most of these problems are 
in fact still open as they turn out to be surprisingly difficult in full
generality. Two outstanding questions among these are the following: let 
$\Om\subseteq\bC$ be an appropriate set of interest and denote by $\pi(\Om)$
the class of all (complex or real) univariate polynomials whose zeros lie in
$\Om$. 

\begin{problem}\label{pb1}
Characterize all linear transformations 
$T:\pi(\Om)\rightarrow \pi(\Om)\cup\{0\}$. 
\end{problem}

Let $\pi_n$ be the vector space (over $\bC$ or $\bR$) of all polynomials
of degree at most $n$ and denote by $\pi_n(\Om)$ the subclass of $\pi(\Om)$
consisting of polynomials of degree at most $n$. The finite degree analog
of Problem~\ref{pb1} is as follows.

\begin{problem}\label{pb2}
Describe all linear operators $T:\pi_n(\Om)\rightarrow \pi(\Om)\cup\{0\}$ for 
$n\in\bN$.
\end{problem}

These problems were stated in precisely this general form in 
\cite{CC1,cso} 
(see also \cite{BBCV}, \cite[pp.~182--183]{RS}) thereby encompassing 
essentially all
similar questions or variations thereof scattered throughout the literature.
Problems~\ref{pb1}--\ref{pb2} originate from the works of 
Laguerre \cite{lag} and P\'olya-Schur \cite{PS} and have been open for all but 
trivial
choices of $\Om$, including such important cases when $\Om=\bR$
or $\Om$ is a half-plane. In this paper we completely solve 
Problems~\ref{pb1}--\ref{pb2} in arguably the most relevant cases, 
namely all closed circular domains (iii)--(v) 
and their boundaries (i)--(ii):
\begin{itemize}
\item[(i)] $\Om$ is a line, 
\item[(ii)] $\Om$ is a circle, 
\item[(iii)] $\Om$ is a closed half-plane, 
\item[(iv)] $\Om$ is a closed disk, 
\item[(v)] $\Om$ is the complement of an open disk.
\end{itemize}

Despite their long history only relatively few results pertaining to
Problems~\ref{pb1}--\ref{pb2} are known. As we note in the following 
(very brief) survey, these deal almost 
exclusively with special
types of linear transformations satisfying the required properties. 

To prove the transcendental characterizations of linear preservers of 
polynomials whose zeros are located on a line or 
in a closed half-plane (Theorems~\ref{tranR} and~\ref{tranC}) we first 
establish a result
on uniform limits on compact sets of bivariate polynomials 
which are non-vanishing whenever  both variables are in the upper half-plane
(Theorem~\ref{genpol}). Entire functions which are 
uniform limits on compact sets 
of sequences of univariate polynomials with only positive zeros were first 
described 
by Laguerre \cite{lag}. In the process he showed that if $Q(z)$ is a real 
polynomial with all negative zeros then 
$T(\pi(\bR))\subseteq \pi(\bR)\cup\{0\}$, 
where $T:\bR[z]\rightarrow \bR[z]$ is the linear operator defined by 
$T(z^k)=Q(k)z^k$, $k\in \bN$. Laguerre also stated without proof 
the correct result  for uniform limits of polynomials with all real zeros.
The class of entire functions thus obtained -- the so-called 
Laguerre-P\'olya class -- was subsequently characterized by P\'olya \cite{pol}.
A more complete investigation of sequences of such polynomials was 
carried out in 
\cite{Lindwart-Polya}. This also led to the description of entire 
functions obtained as uniform limits on compact sets of sequences of 
univariate polynomials all whose zeros lie in a given closed 
half-plane \cite{Le,Ob0,szasz} as well as the description of entire 
functions in two variables obtained 
as limits, uniformly on compact sets, of sequences of bivariate polynomials 
which are non-vanishing when both variables are in a given open half-plane 
\cite{Le}.

The Laguerre-P\'olya class has ever since 
played a significant role in the theory of entire functions 
\cite{CCS,Le}. It was for instance a key ingredient in 
P\'olya and Schur's (transcendental) 
characterization of {\em multiplier sequences of the first kind}  \cite{PS}, 
see Theorem~\ref{ps} below. The latter are linear transformations $T$ on 
$\bR[z]$ that are 
diagonal in the standard monomial basis of $\bR[z]$ and satisfy 
$T(\pi(\bR))\subseteq \pi(\bR)\cup\{0\}$. P\'olya-Schur's seminal paper 
generated a vast literature on this topic and related subjects at the interface
between analysis, operator theory and algebra but a solution to 
Problem~\ref{pb1} in the case $\Om=\bR$ has so far remained elusive 
(cf.~\cite{CC1}). Among the most noticeable progress in this direction we 
should mention Theorem 17 of \cite[Chap.~IX]{Le}, where Levin describes a 
certain class of ``regular" linear operators acting on the closure 
of the set of polynomials in one variable which have all zeros in the closed 
upper half-plane. However, Levin's theorem 
actually uses rather restrictive assumptions and seems in fact to rely on 
additional (albeit not
explicitly stated) non-degeneracy conditions for the transformations involved. 
Indeed, one can easily produce counterexamples to Levin's result by 
considering linear operators such as the ones described in 
Corollary \ref{cor-algstab} (a) of this paper. In \cite{CC2} 
Craven and Csordas established
an analog of the P\'olya-Schur theorem for multiplier sequences in finite 
degree thus solving Problem~\ref{pb2} for $\Om=\bR$ in the special case of
diagonal operators. Unipotent upper triangular linear operators $T$ on 
$\bR[z]$ 
satisfying $T(\pi(\bR))\subseteq \pi(\bR)\cup\{0\}$ were described in 
\cite{CPP}. 
Quite recently, in \cite{BBS} the authors solved Problem~\ref{pb1} for 
$\Om=\bR$ and obtained
multivariate extensions for a large class of linear transformations, namely 
all finite order linear differential operators with polynomial coefficients.
Further partial progress towards a solution to Problem~\ref{pb1} for $\Om=\bR$
is preliminarily reported in \cite{fisk} although the same kind
of remarks as in the case of Levin's theorem apply here. Namely, the
results of {\em op.~cit.~}are valid only in the presence of extra 
non-degeneracy or continuity assumptions for the operators under 
consideration. Various other special cases of Problem~\ref{pb1} for 
$\Om=\bR$ have been 
considered in 
\cite{aleman,Br,brenti, iserles1,iserles2,iserles3}. Finally, we should 
mention that to the best of our 
knowledge  
Problems~\ref{pb1}--\ref{pb2} have so far been widely open in cases 
(ii)--(v).

To begin with, in \S \ref{ss21} and \S \ref{ss31}--\S \ref{ss32} we solve 
Problems~\ref{pb1}--\ref{pb2} for $\Om=\bR$ and 
$\Om=\{z \in \CC: \Im(z) \leq 0\}$ and thus obtain complete algebraic and
transcendental characterizations of linear operators that preserve
hyperbolicity and stability, respectively. In order to deal with 
Problems~\ref{pb1}--\ref{pb2} for {\em all} closed circular domains and their
boundaries we are naturally led to considering a third classification problem,
namely the following more general version of Problem~\ref{pb2}.

\begin{problem}\label{pb3}
Let $n\in\bN$ and $\Omega \subset \CC$. Describe all linear operators 
$$T:\pi_n(\Om)\setminus \pi_{n-1}(\Om) \rightarrow \pi(\Om)\cup\{0\}.$$  
\end{problem}

As we explain in \S \ref{ss22}, Problems~\ref{pb2} and~\ref{pb3} are 
equivalent for closed unbounded sets but for closed bounded sets the latter 
problem is more natural and actually turns out to be a crucial step in 
solving Problems~\ref{pb1}--\ref{pb2} for closed discs. 
In \S \ref{ss22} and \S \ref{ss33} we fully answer Problem~\ref{pb3} for 
any closed circular domain or the boundary of such a domain and as a 
consequence we get
complete solutions to Problems~\ref{pb1}--\ref{pb2} in all cases ((i)--(v))
listed above.

On the one hand, these results accomplish the classification program 
originating from the works of Laguerre and P\'olya-Schur that we briefly 
outlined 
in this introduction. On the other hand, they seem to have numerous 
applications ranging from entire function theory and operator theory to 
real algebraic geometry, matrix theory and combinatorics. Some of these 
will make the object of forthcoming publications.  
The paper concludes with several remarks
on related open problems and potential further developments (\S \ref{s4}).

\section{Main Results}\label{s2}

\subsection{Hyperbolicity and Stability Preservers}\label{ss21}

To formulate the complete answers to Problems~\ref{pb1}--\ref{pb2} for $\RR$ 
and 
the half-plane $\{z \in \CC: \Im(z) \leq 0\}$  we need to introduce some 
notation.
As in \cite{BBS} -- and following the commonly used terminology in 
e.g.~the theory of partial differential equations \cite{ABG} -- we call a 
non-zero univariate polynomial with real 
coefficients {\em hyperbolic} if all its zeros are real. Such a polynomial
is said to be {\em strictly hyperbolic} if in addition 
all its zeros are distinct. 
A univariate polynomial 
$f(z)$ with complex coefficients is called {\em stable} if 
$f(z) \neq 0$ for all $z \in \CC$ with $\Im(z) >0$ and it is called 
{\em strictly stable} if 
$f(z) \neq 0$ for all $z \in \CC$ with $\Im(z) \geq 0$.  Hence a univariate 
polynomial with real coefficients is stable if and only if it is hyperbolic. 
These classical concepts have several natural extensions to multivariate 
polynomials, the most general notion being as follows.

\begin{definition}\label{def1}
A polynomial $f(z_1, \ldots, z_n)\in\bC[z_1,\ldots,z_n]$ is 
{\em stable} 
if  $f(z_1, \ldots, z_n) \neq 0$ for all $n$-tuples 
$(z_1, \ldots, z_n) \in \CC^n$ with $\Im(z_j) >0$, $1\le j\le n$.  
If in addition $f$ has real coefficients it will be referred  to as 
{\em real stable}. The 
sets of stable and real stable polynomials in $n$ variables are denoted 
by $\HH_n(\CC)$
 and 
$\HH_n(\RR)$, respectively. 
\end{definition}

Note that $f$ is stable (respectively, real stable) if and only if for all 
$\alpha \in \RR^n$ and $v \in \RR_+^n$ the 
univariate polynomial $f(\alpha + vt)\in\bC[t]$ is stable 
(respectively, hyperbolic). The connection between real stability and 
(G\aa rding) hyperbolicity
for multivariate homogeneous polynomials is explained in 
e.g.~\cite[Proposition 1]{BBS}.

\begin{notation}\label{not-cplx}
Henceforth it is understood that if $T$ is a linear operator on
some (real) linear subspace $V\subseteq \bR[z_1,\ldots,z_n]$ then $T$ extends 
in an obvious fashion to a linear operator -- denoted again by $T$ -- on the 
complexification $V\oplus iV$ of $V$.
\end{notation}

\begin{definition}\label{def2}
A linear operator $T$ defined on a linear subspace $V$ of 
$\bC[z_1,\ldots,z_n]$ (respectively,
$\bR[z_1,\ldots,z_n]$) is called {\em stability preserving} 
(respectively, {\em real stability preserving}) on a given subset 
$M\subseteq V$ if
$$T(\HH_n(\CC)\cap M)\subseteq \HH_n(\CC)\cup\{0\}\text{ (respectively,
$T(\HH_n(\bR)\cap M)\subseteq \HH_n(\bR)\cup\{0\}$)}.$$ 
A real stability preserver 
in the univariate case will alternatively be referred to as a 
{\em hyperbolicity preserver}. For $m\in\bN$ let 
$\RR_m[z] = \{ f \in \RR[z] : \deg(f) \leq m\}$ and 
$\CC_m[z]=\RR_m[z]\oplus i\RR_m[z]=\{ f \in \CC[z] : \deg(f) \leq m\}$. 
If $T$ is a stability (respectively, hyperbolicity) preserving operator 
on $\bC_m[z]$ (respectively, $\bR_m[z]$)
we will also say that $T$ {\em preserves stability} (respectively, 
{\em hyperbolicity}) {\em up to degree} $m$.
\end{definition}

P\'olya-Schur's characterization of multiplier sequences of the first kind
that we already alluded to in the introduction is given in the following 
theorem \cite{CC1,Le,PS}.

\begin{theorem}[P\'olya-Schur theorem]\label{ps}
Let $\lambda : \NN \rightarrow \RR$ be a sequence of real numbers and 
$T: \RR[z] \rightarrow \RR[z]$ be the corresponding (diagonal) linear operator
given by $T(z^n)=\la(n)z^n$, $n\in\bN$.
 Define $\Phi(z)$ to be the formal power series 
$$
\Phi(z) = \sum_{k=0}^\infty \frac{\lambda(k)}{k!}z^k.
$$
The following assertions are equivalent:
\begin{itemize}
\item[(i)] $\lambda$ is a multiplier sequence,  
\item[(ii)] $\Phi(z)$ defines an entire 
function which is the limit, uniformly on compact sets, of 
polynomials with only real zeros of the same sign, 
\item[(iii)] Either $\Phi(z)$ or $\Phi(-z)$ is an entire function 
that can be written as
$$C z^n e^{az} \prod_{k=1}^\infty (1+ \alpha_k z),$$ 
where $n \in \NN$, $C \in \RR$, $a,\alpha_k \geq 0$ for all $k \in \NN$ and 
$\sum_{k=1}^\infty \alpha_k < \infty$, 
\item[(iv)] For all non-negative integers $n$ the  polynomial 
$T[(z+1)^n]$ is hyperbolic with all zeros of the same sign. 
\end{itemize}
\end{theorem} 

As noted in e.g.~\cite[Theorem 3.3]{CC1}, parts (ii)-(iii) in P\'olya-Schur's 
theorem give a 
``transcendental'' description of multiplier sequences while part (iv) 
provides an ``algebraic'' characterization. We emphasize right away the fact
that our main results actually yield 
algebraic and transcendental characterizations of {\em all}
hyperbolicity and stability preservers, respectively, and are therefore
natural generalizations of Theorem~\ref{ps}. Moreover, they also display
an intimate connection between Problem~\ref{pb1} and its finite degree analog
(Problem~\ref{pb2}) in the case of (real) stability preservers. 

\begin{notation}\label{not1}  
Given a linear operator  
$T$ on $\CC[z]$ we extend it to a linear operator -- denoted
again by $T$ -- on the space $\bC[z,w]$ of polynomials in 
the variables $z,w$ by setting
$T(z^kw^\ell) = T(z^k)w^\ell$
for all $k,\ell \in \NN$.
\end{notation}

\begin{definition}\label{def3}
Let $\alpha_1 \leq \alpha_2 \leq \cdots \leq \alpha_n$ and 
$\beta_1 \leq \beta_2 \leq \cdots \leq \beta_m$ be the zeros of two 
hyperbolic polynomials $f,g \in \HH_1(\RR)$. We say that these zeros 
{\em interlace} if they can 
be ordered so that either   
$\alpha_1 \leq \beta_1 \leq \alpha_2 \leq \beta_2 \leq 
\cdots$ or $\beta_1 \leq \alpha_1 \leq \beta_2 \leq \alpha_2 \leq 
\cdots$. Note that in this case one has $|m-n|\le 1$. By convention, 
the zeros of any two polynomials 
of degree $0$ or $1$ interlace. 
\end{definition} 

Our first theorem characterizes linear operators preserving hyperbolicity up 
to some fixed degree $n$. 

\begin{theorem}\label{finitehyp}
Let $n \in \NN$ and let $T : \RR_n[z] \rightarrow \RR[z]$ be a linear 
operator. Then $T$ preserves hyperbolicity 
if and only if either
\begin{itemize}
\item[(a)] $T$ has range of dimension at most two and is of the form 
$$
T(f) = \alpha(f)P + \beta(f)Q,\quad f\in\bR_n[z],
$$
where $\alpha, \beta: \RR_n[z] \rightarrow \RR$ are linear functionals and 
$P, Q\in\HH_1(\bR)$ have interlacing zeros, or 
\item[(b)] $T[(z+w)^n]\in \HH_2(\RR)$, or 
\item[(c)] $T[(z-w)^n]\in \HH_2(\RR)$. 
\end{itemize}
\end{theorem} 

Real stable polynomials in two variables have recently been characterized by 
the authors 
\cite{BBS} as the polynomials $f(z,w)\in\bR[z,w]$ that can be expressed as
\begin{equation}\label{lax}
f(z,w)=\pm \det(zA+wB+C),
\end{equation}
where $A$ and $B$ are positive semi-definite matrices and $C$ is a symmetric 
matrix.  Hence 
(b) and (c) in Theorem~\ref{finitehyp} can be reformulated as 
$$
T[(z+w)^n]= \pm \det(zA \pm wB+C),
$$
where $A$ and $B$ are positive semi-definite matrices and $C$ is a symmetric 
matrix. 

We will also need to deal with the case when we allow complex coefficients. 
\begin{theorem}\label{finitehypC}
Let $n \in \NN$ and let $T : \CC_n[z] \rightarrow \CC[z]$ be a linear 
operator. Then 
$T: \pi_n(\RR) \rightarrow \pi(\RR)$ 
if and only if either 
\begin{itemize}
\item[(a)] $T$ has range of dimension at most one  and is of the form 
$$
T(f) = \alpha(f)P, \quad f\in\bC_n[z],
$$
where $\alpha:\CC_n[z] \rightarrow \CC$ is a linear functional and 
$P\in\HH_1(\bR)$, or 
\item[(b)] $T$ has range of dimension at most two and is of the form 
$$
T(f) = \eta\alpha(f)P + \eta\beta(f)Q,\quad f\in\bC_n[z],
$$
where $\eta \in \CC$,  $\alpha, \beta: \CC_n[z] \rightarrow \CC$  are linear 
functionals such that $\al(\bR_n[z])\subseteq \bR$, $\be(\bR_n[z])\subseteq 
\bR$, and $P,Q\in\HH_1(\bR)$ have interlacing zeros, or 
\item[(c)] There exists $\eta\in\bC$ such that  
$\eta T[(z+w)^n] \in \HH_2(\RR)$, or 
\item[(d)] There exists $\eta\in\bC$ such that  
$\eta T[(z-w)^n]\in \HH_2(\RR)$. 
\end{itemize}
\end{theorem} 

The corresponding theorem for stability preservers up to some fixed degree 
$n$ reads as follows. 

\begin{theorem}\label{finitestab}
Let $n \in \NN$ and let $T : \CC_n[z] \rightarrow \CC[z]$ be a linear 
operator. Then $T$ preserves stability 
if and only if either
\begin{itemize}
\item[(a)] $T$ has range of dimension at most one and is of the form 
$$
T(f) = \alpha(f)P,\quad f\in\bC_n[z],
$$
where $\alpha:\CC_n[z] \rightarrow \CC$ is a linear functional and 
$P\in\HH_1(\bC)$, or 
\item[(b)] $T[(z+w)^n]\in \HH_2(\CC)$. 
\end{itemize}
\end{theorem}

From Theorems~\ref{finitehyp} and~\ref{finitestab} we deduce the following
algebraic characterizations of hyperbolicity and stability preservers,
respectively.

\begin{corollary}[Algebraic Characterization of Hyperbolicity 
Preservers]\label{cor-alghyp}
A linear operator $T:\bR[z]\rightarrow \bR[z]$ preserves hyperbolicity
if and only if either
\begin{itemize}
\item[(a)] $T$ has range of dimension at most two and is of the form 
$$
T(f) = \alpha(f)P + \beta(f)Q,\quad f\in\bR[z],
$$
where $\alpha,\beta:\bR[z]\rightarrow \bR$ are linear functionals and 
$P,Q\in\HH_1(\bR)$ have interlacing zeros, or 
\item[(b)]
$T[(z + w)^n]\in \HH_2(\bR)\cup \{0\}$ for all $n\in \bN$, or
\item[(c)]
$T[(z -w)^n]\in \HH_2(\bR)\cup \{0\}$ for all $n\in \bN$.
\end{itemize}
\end{corollary}

\begin{corollary}[Algebraic Characterization of Stability 
Preservers]\label{cor-algstab}
A linear operator $T:\bC[z]\rightarrow \bC[z]$ preserves stability 
if and only if either
\begin{itemize}
\item[(a)] $T$ has range of dimension at most one and is of the form 
$$
T(f) = \alpha(f)P,\quad f\in\bC[z],
$$
where $\alpha:\bC[z]\rightarrow \bC$ is a linear functional and 
$P\in\HH_1(\bC)$, or 
\item[(b)] $T[(z+w)^n]\in \HH_2(\bC)\cup \{0\}$ for all $n\in \bN$.
\end{itemize}
\end{corollary} 

\begin{notation}\label{not2}
To any linear operator $T:  \CC[z] \rightarrow \CC[z]$ we associate a formal 
power series in $w$ 
with polynomial coefficients in $z$
$$
G_T(z,w)= \sum_{n=0}^\infty \frac{(-1)^nT(z^n)}{n!}w^n\in\bC[z][[w]]. 
$$
Let $n$ be a positive integer and denote by $\HC_n(\CC)$ and $\HC_n(\RR)$, 
respectively, the set of entire functions in $n$ variables that are limits, 
uniformly on compact sets, 
of  polynomials in $\HH_n(\CC)$ and $\HH_n(\RR)$, respectively. Hence in our 
notation $\HC_1(\RR)$ is the Laguerre-P\'olya class of  entire functions, 
sometimes denoted by $\mathcal{L}$-$\mathcal{P}$ in the literature. For a 
description of $\HC_1(\CC)$ and $\HC_2(\CC)$ we refer 
to \cite[Chap. IX]{Le}.
\end{notation} 

\begin{remark}
As noted in \cite{BBS}, any linear operator $T$ on $\bC[z]$ may be uniquely 
represented as a formal linear differential operator with polynomials 
coefficients, i.e., $T=\sum_{k=0}^{\infty}Q_k(z)\frac{d^k}{d z^k}$, where
$Q_k\in \bC[z]$, $k\ge 0$. In \cite{BBS} we used the (formal) 
symbol of $T$, i.e., 
$F_T(z,w):=\sum_{k=0}^{\infty}Q_k(z)w^k\in\bC[z][[w]]$. One can easily check
that the
``modified symbol'' $G_T(z,w)$ introduced in Notation~\ref{not2} satisfies
$G_T(z,w)e^{zw}=F_T(z,-w)$.
\end{remark}

\begin{theorem}[Transcendental Characterization of Hyperbolicity 
Preservers]\label{tranR}
A linear operator $T: \RR[z] \rightarrow \RR[z]$ preserves hyperbolicity if 
and only if either
\begin{itemize}
\item[(a)] $T$ has range of dimension at most two and is of the form 
$$
T(f) = \alpha(f)P + \beta(f)Q,\quad f\in\bR[z],
$$
where $\alpha,\beta:\bR[z]\rightarrow \bR$ are linear functionals and 
$P,Q\in\HH_1(\bR)$ have interlacing zeros, or 
\item[(b)] $G_T(z,w) \in \HC_2(\RR)$, or 
\item[(c)]  $G_T(z,-w)\in \HC_2(\RR)$. 
\end{itemize}
\end{theorem}

\begin{theorem}[Transcendental Characterization of Stability 
Preservers]\label{tranC}
A linear operator $T: \CC[z] \rightarrow \CC[z]$ preserves stability if and 
only if either
\begin{itemize}
\item[(a)] $T$ has range of dimension at most one and is of the form 
$$
T(f) = \alpha(f)P, \quad f\in\bC[z],
$$
where $\alpha:\CC[z] \rightarrow \CC$ is a linear functional and 
$P\in\HH_1(\bC)$, or 
\item[(b)] $G_T(z,w) \in \HC_2(\CC)$.
\end{itemize}
\end{theorem}

\subsection{Preservers of Polynomials with Zeros in a Closed Circular Domain 
or Its Boundary}\label{ss22}

Recall that a M\"obius transformation is a bijective conformal map of the 
extended complex plane, i.e., 
a map $\Phi : \CC\cup \{\infty\} \rightarrow \CC \cup \{\infty\}$ given by  
\begin{equation}\label{mobius}
\Phi(z)= \frac {az+b}{cz+d}, \quad a,b,c,d \in \CC,\, ad-bc \neq 0. 
\end{equation}
The inverse of $\Phi$ is then given by 
$$ 
 \Phi^{-1}(z)= \frac{dz-b}{-cz+a}. 
$$

\begin{definition}\label{def-circ}
Let $H=\{z \in \CC : \Im(z) >0\}$ and 
$\overline{H}=\{z \in \CC : \Im(z) \ge 0\}$. An {\em open circular domain} is 
the image 
of $H$ under a  M\"obius transformation, i.e., an open disk, the (open) 
complement of a closed disk or an 
open affine half-plane. A {\em closed circular domain} is the image 
of $\overline{H}$ under a M\"obius transformation, that is, a closed disk, 
the (closed) complement of an open disk or a
closed affine half-plane. 
\end{definition}

For technical reasons we will henceforth assume the following:
\begin{equation}\label{assum}
\begin{split}
&\text{{\em If $C$ is a half-plane then the corresponding M\"obius 
transformation}}\\ 
&\text{{\em $\Phi : C \rightarrow H$  is a translation composed with a 
rotation, i.e., $c=0$ in \eqref{mobius}.}}
\end{split}
\end{equation}

Let us also extend Definition~\ref{def1} to arbitrary sets 
$\Omega\subseteq \CC$. 

\begin{definition}\label{def-stab-om}
A polynomial $f \in \CC[z_1,\ldots, z_n]$ is called 
$\Omega$-{\em stable} if $f(\zeta_1, \ldots, \zeta_n)\neq 0$ whenever 
$\zeta_j \in \Omega$ for all 
$1 \leq j \leq n$. 
\end{definition} 

\begin{remark}\label{r-h-stab}
Note that an $H$-stable polynomial is precisely a stable polynomial in the
sense of Definition~\ref{def1}.
\end{remark}

A fundamental discrepancy between $\pi_n(C_1)$ and $\pi_n(C_2)$, where $C_1$ 
is closed and unbounded 
and $C_2$ is closed and bounded, is that 
$\pi_n(C_1) \setminus \pi_{n-1}(C_1)$ is dense in $\pi_n(C_1)$ while  
$\pi_n(C_2) \setminus \pi_{n-1}(C_2)$ is not dense in $\pi_n(C_2)$ since 
constant non-zero polynomials do not belong to
$\overline{\pi_n(C_2) \setminus \pi_{n-1}(C_2)}$.  In order for a linear 
transformation $T: \CC_n[z] \rightarrow \CC[z]$ to map  
$\pi_n(C_1) \cup \{0\}$  into  $\pi(C_1) \cup\{0\}$ it is therefore enough 
(by Hurwitz' theorem) for 
$T$ to map $\pi_n(C_1) \setminus \pi_{n-1}(C_1)$ into  $\pi(C_1) \cup\{0\}$. 
However, this is not the case for $C_2$.  Indeed, take for instance $C_2$ to 
be the closed unit disk and let $T_n : \CC_n[z] \rightarrow \CC_{n-1}[z]$ be 
defined by 
$$
T(z^k) = (n-k)z^k + kz^{k-1}, \quad 0\le k\le n.
$$ 
An application of Theorem~\ref{circularregion} below shows that 
$T: \pi_n(C_2) \setminus \pi_{n-1}(C_2) \rightarrow \pi(C_2) \cup\{0\}$ but 
$T(z^{n-1})=z^{n-2}(z+n-1)$, which does not belong to $\pi(C_2)$ for 
$n \geq 3$. Therefore we need to solve a more general version of 
Problem~\ref{pb2}, namely Problem~\ref{pb3} in \S \ref{s1}. This is done in 
the next two theorems for all closed circular domains and their boundaries.

\begin{notation}\label{not-sets}
Given $\Omega \subseteq \CC$  we denote its
complement $\bC\setminus \Omega$ by $\Omega'$, its boundary 
$\overline{\Om}\setminus\mathring{\Om}$ by $\partial\Om$, and we let $\Om^r$ 
be the interior of
the complement of $\Om$, that is, $\Om^r=\mathring{\Om'}$. 
\end{notation}

\begin{theorem}\label{circularregion}
Let $n \in \NN$ and $T : \CC_n[z] \rightarrow \CC[z]$ be a linear operator. Let
further $C$ be an open circular domain given by $C=\Phi^{-1}(H)$,
where $\Phi$ 
is a M\"obius transformation as in \eqref{mobius}. 
Then $T : \pi_n(C')\setminus \pi_{n-1}(C') \rightarrow \pi(C')\cup\{0\}$ if 
and only if either
\begin{itemize}
\item[(a)] $T$ has range of dimension at most one and is of the form 
$$
T(f) = \alpha(f)P,\quad f\in\CC_n[z],
$$
where $\alpha: \CC_n[z] \rightarrow \CC$ is a linear functional and 
$P \in \pi(C')$, or 
\item[(b)] The polynomial 
$$
T\left[\big((az+b)(cw+d)+(aw+b)(cz+d)\big)^n\right] 
$$
is $C$-stable.
\end{itemize}
\end{theorem}
\begin{remark}
In the case when $C'$ is the closed unit disk then for 
\begin{equation}\label{spec-fi}
\Phi(z)=\frac{(i/2)(z+i)}{z-i}
\end{equation}
the polynomial in (b) of Theorem~\ref{circularregion} reduces to 
$i^nT[(1+zw)^n]$. 
\end{remark}

\begin{notation}\label{not-fi}
Given a M\"obius transformation $\Phi$ as in \eqref{mobius} and $n\in\bN$ 
we define an invertible linear transformation  
$\phi_n : \CC_n[z] \rightarrow \CC_n[z]$  by 
$\phi_n(f)(z) = (cz+d)^nf(\Phi(z))$.
\end{notation}
 
\begin{theorem}\label{circle}
Let $n \in \NN$ and $T : \CC_n[z] \rightarrow \CC[z]$ be a linear operator 
and let further $C=\Phi^{-1}(H)$ be an unbounded open circular domain, where 
$\Phi$ 
is a M\"obius transformation as in \eqref{mobius}. 
Then 
$T : \pi_n(\partial C)\setminus \pi_{n-1}(\partial C) \rightarrow 
\pi(\partial C)\cup\{0\}$ if and only if either
\begin{itemize}
\item[(a)]  $T$ has range of dimension at most one and is of the form 
$$
T(f) = \alpha(f)P, \quad f\in\CC_n[z],
$$
where $\alpha: \CC_n[z] \rightarrow \CC$ is a linear functional and 
$P \in \pi(\partial C)$, or 
\item[(b)] $T$ has range of dimension  two and the linear operator given 
by $S=\phi^{-1}_mT\phi_n$ is a stability preserver as in (b) of 
Theorem~\ref{finitehypC}, where $m = \max\{\deg T(f): f \in \CC_n[z]\}$, or
\item[(c)] The polynomial  
$$
T\left[\big((az+b)(cw+d)+(aw+b)(cz+d)\big)^n\right] 
$$
is both $C$-stable and $C^r$-stable,  or 
\item[(d)] The polynomial 
$$
T\left[\big((az+b)(cw+d)-(aw+b)(cz+d)\big)^n\right]
$$
is both $C$-stable and $C^r$-stable.
\end{itemize}
\end{theorem}

\begin{remark}
If $\partial C$ is the unit circle then for $\Phi$ as in~\eqref{spec-fi} the 
polynomials in (c) and (d) of Theorem~\ref{circle} simply become 
$i^nT[(1+zw)^n]$ and $T[(z-w)^n]$, respectively. 
\end{remark}

For the sake of completeness we also formulate analogs of 
Corollaries~\ref{cor-alghyp} and~\ref{cor-algstab} providing algebraic 
characterizations in the case of closed circular domains and their boundaries. 

\begin{corollary}[Algebraic Characterization: Closed Circular Domain 
Case]\label{alg-circularregion}
Let $T : \CC[z] \rightarrow \CC[z]$ be a linear operator and let
$C\subset \CC$ be an open circular domain given by $C=\Phi^{-1}(H)$,
where  $\Phi$ 
is a M\"obius transformation as in \eqref{mobius}. 
Then $T : \pi(C') \rightarrow \pi(C')\cup\{0\}$ if and only if either
\begin{itemize}
\item[(a)] $T$ has range of dimension at most one and is of the form 
$$
T(f) = \alpha(f)P,\quad f\in \bC[z],
$$
where $\alpha: \CC[z] \rightarrow \CC$ is a linear functional and 
$P \in \pi(C')$, or 
\item[(b)] For all $n \in \NN$ the polynomial 
$$
T\left[\big((az+b)(cw+d)+(aw+b)(cz+d)\big)^n\right] 
$$
is $C$-stable.
\end{itemize}
\end{corollary}
 
\begin{corollary}[Algebraic Characterization: Circle and Line 
Case]\label{cor-circle}
Let  $T : \CC[z] \rightarrow \CC[z]$ be a linear operator and 
$C=\Phi^{-1}(H)$ be an unbounded open circular domain,  where $\Phi$ 
is a M\"obius transformation as in \eqref{mobius}. 
Then $T : \pi(\partial C) \rightarrow \pi(\partial C)\cup\{0\}$ if and only 
if either
\begin{itemize}
\item[(a)]  $T$ has range of dimension at most one and is of the form 
$$
T(f) = \alpha(f)P,\quad f\in \bC[z],
$$
where $\alpha: \CC[z] \rightarrow \CC$ is a linear functional and 
$P\in\pi(\partial C)$, or 
\item[(b)] $T$ has range of dimension  two and  for all $n \in \NN$ the 
linear operator given by $S_n=\phi^{-1}_{m(n)}T\phi_n$ is a stability 
preserver as in (b) of Theorem~\ref{finitehypC}, where 
$m(n) = \max\{\deg T(f): f \in \CC_n[z]\}$, or  
\item[(c)] For all $n \in \NN$ the polynomial  
$$
T\left[\big((az+b)(cw+d)+(aw+b)(cz+d)\big)^n\right] 
$$
is both $C$-stable and $C^r$-stable,  or 
\item[(d)] For all $n \in \NN$ the polynomial 
$$
T\left[\big((az+b)(cw+d)-(aw+b)(cz+d)\big)^n\right]
$$
is both $C$-stable and $C^r$-stable.
\end{itemize}
\end{corollary}

Similarly, we may characterize all linear maps that take polynomials with 
zeros in one closed 
circular domain $\Om_1$ to polynomials with zeros in another closed circular 
domain $\Om_2$, or 
the boundary of one  
circular domain $\Om_1$ to the boundary of another circular domain $\Om_2$. 
However, this only amounts to composing with linear operators of the type 
defined
in Notation~\ref{not-fi}, namely 
$\phi_n : \CC_n[z] \rightarrow \CC_n[z]$, where 
$
\phi_n(f)(z) = (cz+d)^n f(\Phi(z))
$ 
and $\Phi$ is an appropriate M\"obius transformation of the form 
\eqref{mobius}. 

\section{Proofs of the Main Results}\label{s3}

\subsection{Hyperbolic and Stable Polynomials}\label{ss31}

If the zeros of two hyperbolic polynomials $f,g\in\HH_1(\bR)$ interlace then 
the {\em Wronskian} 
$W[f,g]:=f'g-fg'$ is either non-negative or non-positive on the whole real 
axis $\RR$. 

\begin{definition}\label{def-pp}
Given $f,g \in \HH_1(\bR)$ we say that $f$ and $g$ 
are in {\em proper position}, denoted $f \ll g$, if the zeros of $f$ and $g$ 
interlace and $W[f,g] \leq 0$. 
\end{definition}

For technical reasons we also say that the zeros of the polynomial $0$ 
interlace the zeros of any (non-zero) hyperbolic polynomial and write 
$0 \ll f$ and $f \ll 0$. Note 
that if $f \ll g$ and $g \ll f$ then $f$ and $g$ must be constant multiples 
of each other, that is, $W[f,g]\equiv 0$.

The following theorem is a version of the classical Hermite-Biehler 
theorem \cite{RS}.

\begin{theorem}[Hermite-Biehler theorem]
Let $h := f +ig \in \CC[z]$, where $f,g \in \RR[z]$. Then 
$h \in \HH_1(\CC)$ if and only if $f,g\in\HH_1(\bR)$ and $g \ll f$. 
Moreover, $h$ is strictly stable  if and only if 
$f$ and $g$ are strictly hyperbolic polynomials with no common zeros and 
$g \ll f$.
\end{theorem}  

The next theorem is often attributed to Obreschkoff \cite{obreschkoff}.
 
\begin{theorem}[Obreschkoff theorem]
Let $f,g \in \RR[z]$. Then 
$\alpha f + \beta g \in \HH_1(\RR)\cup \{0\}$ for all $\alpha, \beta \in \RR$ 
if and only if either 
$f \ll g$, $g \ll f$, or  $f=g\equiv 0$.  
Moreover,  $\alpha f + \beta g$ is strictly hyperbolic for all 
$\alpha, \beta \in \RR$ with 
$\alpha^2 + \beta^2 \neq 0$ if and only if  $f$ and $g$ are strictly 
hyperbolic polynomials with no common zeros and either $f \ll g$ or $g \ll f$. 
\end{theorem}

\begin{remark}\label{obremark}
Note that if $T : \pi_n(\RR) \rightarrow \pi(\RR)$ is an $\RR$-linear 
operator then by Obreschkoff's theorem  $T$ also preserves 
interlacing in the following manner: if $f$ and $g$ are hyperbolic 
polynomials of degree 
at most $n$ whose zeros interlace then the zeros of $T(f)$ and $T(g)$ 
interlace 
provided that $T(f)T(g) \neq 0$. 
\end{remark} 

\begin{lemma}\label{degenerate}
Let $n \in \NN$. 
Suppose that  $T : \RR_{n+1}[z] \rightarrow \RR[z]$ preserves hyperbolicity 
and 
that $f \in \RR[z]$ is a strictly hyperbolic polynomial of degree $n$ or 
$n+1$ for which $T(f)=0$. Then 
$T(g)$ is hyperbolic for all $g \in \RR[z]$ with $\deg g \leq n+1$. 

Let $T : \CC_n[z] \rightarrow \CC[z]$ be a stability preserver and suppose 
that $f \in \CC[z]$ is a strictly stable polynomial of degree $n$ for which 
$T(f)=0$. Then 
$T(g)$ is stable for all $g \in \CC[z]$ with $\deg g \leq n$. 
\end{lemma}

\begin{proof}
Let $f$ be a strictly hyperbolic polynomial of degree $n$ or $n+1$ for which 
$T(f)=0$ and let $g \in \RR[z]$ be a 
polynomial with $\deg g \leq n+1$. From Hurwitz' theorem it follows that for 
$\epsilon \in \RR$ with $|\epsilon|$ small enough 
the polynomial $f+\epsilon g$ is strictly hyperbolic. Since 
$\deg (f+\eps g)\le n+1$ and $T$ preserves hyperbolicity up to degree $n+1$ 
we get that $T(g)=\epsilon^{-1} T(f+\eps g)$ is 
hyperbolic. 

Suppose that  $f \in \CC[z]$ is a strictly stable polynomial of degree $n$ 
such that 
$T(f)=0$ and that the degree of $g \in \CC[z]$ does not exceed $n$. By 
Hurwitz' theorem 
$f+\epsilon g$ is strictly stable for all sufficiently small $|\epsilon|$. 
Since $\deg (f+\eps g)\le n$ and $T$ preserves stability up to degree $n$ it 
follows that $T(g)=\epsilon^{-1} T(f+\epsilon g)$ is stable. 
\end{proof}

\begin{lemma}\label{spaces}
Suppose that $V \subseteq \RR[z]$ is an $\RR$-linear space whose every 
non-zero element is hyperbolic. Then $\dim V \leq 2$. 

Suppose that $V \subseteq \CC[z]$ is a $\CC$-linear space whose every 
non-zero element is stable. Then $\dim V \leq 1$. 
\end{lemma}

\begin{proof}
We first deal with the real case. Suppose that there are three linearly 
independent polynomials $f_1,f_2$  and $f_3$ 
in $V$. By Obreschkoff's theorem the zeros of these polynomials mutually 
interlace. Wlog we may assume that 
$f_1 \ll f_2$ and $f_1 \gg f_3$. Consider the line  segment 
$\ell_\theta =\theta f_3 + (1-\theta) f_2$, 
$0 \leq \theta \leq 1$. Since $f_1 \ll \ell_0$ and $f_1 \gg \ell_1$ by 
Hurwitz' theorem  
there is a real number $\eta$ between $0$ and $1$ such that 
$f_1 \ll \ell_\eta$ and $f_1 \gg \ell_\eta$. This 
means that $f_1$ and $\ell_\eta$ are constant multiples of each other 
contrary to the assumption that 
$f_1,f_2$  and $f_3$ are linearly independent. 

For the complex case 
let $V_R= \{ p : p+iq \in V\text{ with } p,q\in \RR[z] \}$ be the 
``real component'' of $V$. By the Hermite-Biehler theorem 
all polynomials in $V_R$ are hyperbolic,  so by the above we have 
$\dim_{\bR} V_R \leq 2$. Clearly, 
$V$ is the complex span of $V_R$. If $\dim_{\bR} V_R \leq 1$ we are done 
so we may assume that 
$\{p,q\}$ is a basis for $V_R$ with $f:=p+iq \in V$. By definition 
$W[p,q] \geq 0$ on the whole of $\RR$ and 
the Wronskian is not identically zero. Assume now that $g$ is another 
polynomial in $V$. Then 
$$
g= ap+bq+i(cp+dq) 
$$
for some $a,b,c,d \in \RR$. We have to show that $g$ is a (complex) constant 
multiple of $f$. Since $g \in V$ we have  
$$
W[ap+bq,cp+dq] = (ad-bc)W[p,q] \geq 0,
$$
so that $ad-bc \geq 0$. Now by linearity we have that 
$$g+(u+iv)f =(a+u)p+(b-v)q+i( (c+v)p+(d+u)q)    \in V
$$ 
for all $u,v \in \RR$ which, as above, gives 
$$
H(u,v):=(a+u)(d+u)-(b-v)(c+v) \geq 0 
$$
for all $u,v \in \RR$. But 
$$
4H(u,v)=(2u+a+d)^2+(2v+c-b)^2-(a-d)^2-(b+c)^2,
$$
so $H(u,v) \geq 0$ for all $u,v \in \RR$ if and only if $a=d$ and $b=-c$. 
This gives  
$$
g=ap-cq+i(cp+aq)=(a+ic)(p+iq)=(a+ic)f,
$$
as was to be shown.
 \end{proof}

\begin{notation}\label{not4}
Let $\HH_1^-(\CC)=\{f \in \CC[z]:f(z) \neq 0\text{ if } \Im(z)<0\}$. By the 
Hermite-Biehler theorem and Definition~\ref{def-pp} if $f,g \in \RR[z]$ then 
$f+ig \in \HH_1^-(\CC)$ if and only if $f,g\in\HH_1(\RR)$ and $f \ll g$. 
Given a linear subspace $V$ of $\CC[z]$ and $M\subseteq V$ 
we say that a linear operator $T$  
on $V$ is {\em stability reversing} on $M$ if 
$T(\HH_1(\CC)\cap M)\subseteq \HH_1^-(\CC) \cup \{0\}$. Note that if 
$S : \CC[z] \rightarrow \CC[z]$ is the 
linear involution defined by $S(f)(z)= f(-z)$ then $T$ is stability reversing 
if and only if $S\circ T$ is stability  preserving. 
\end{notation}

\begin{lemma}\label{trichotomy}
Suppose that  $T : \RR_n[z] \rightarrow \RR[z]$  maps all hyperbolic 
polynomials of degree at most $n$ 
to hyperbolic polynomials. Then 
$T$ is either stability preserving or stability reversing or the range of 
$T$ has dimension at most two. In the latter case 
$T$ is given by 
\begin{equation}\label{form2}
T(f) = \alpha(f)P + \beta(f)Q,\quad f\in\bR_n[z],
\end{equation}
where $P,Q$ are hyperbolic polynomials whose zeros interlace and 
$\alpha, \beta $ are real-valued linear functionals on $\RR_n[z]$.  
\end{lemma}

\begin{proof} 
The lemma is obvious for $n=0$ so we may and do assume that $n$ is a positive 
integer. 
By Remark \ref{obremark} and the Hermite-Biehler theorem we know that $T$ 
maps all
 stable polynomials of degree $n$ into 
the set $\HH_1(\CC)\cup \HH_1^-(\CC) \cup\{0\}$. We now distinguish two 
cases. Suppose first that there are two strictly 
stable polynomials $f,g$ of degree $n$ such that $T(f) \in \HH_1(\CC)$ and 
$T(g) \in \HH_1^-(\CC)$. 

\medskip

{\em Claim.} If the above conditions are satisfied then the kernel of $T$ 
must contain a strictly hyperbolic polynomial of degree at least $n-1$.

\medskip

From Lemma \ref{degenerate} and the Claim we deduce that 
$T : \RR_n[z] \rightarrow \pi(\RR)\cup\{0\}$. Hence 
all non-zero polynomials in the image of $T$ are hyperbolic, which by 
Lemma \ref{spaces} 
gives that $\dim T(\RR_n[z]) \leq 2$ and thus $T$ must be of the 
form~\eqref{form2}. 

\medskip

{\em Proof of Claim.} Suppose that $f_1,f_2$ are two strictly stable 
polynomials of degree $n$ 
for which 
$T(f_1) \in \HH_1(\CC)$ and $T(f_2) \in \HH_1^-(\CC)$, respectively. By a 
homotopy argument, invoking again 
Hurwitz' theorem, there is a strictly stable polynomial $h$ of degree $n$ for 
which 
$T(h) \in \HH_1(\CC)\cap \HH_1^-(\CC) \cup\{0\}$. Writing $h$ as $h=p+iq$, 
where $p$ and $q$ are 
strictly hyperbolic polynomials (by the Hermite-Biehler theorem) gives that 
$T(p)$ and $T(q)$ are constant multiples of each other. Suppose that 
$\deg p =n$.  Then $\deg q \geq n-1$ since the zeros of $p$ and $q$ 
interlace. If  $T(q)=0$ 
the claim is obviously true so suppose that $T(q) \neq 0$. Then 
$T(p)= \lambda T(q)$ for some $\lambda \in \RR$. 
By the Obreschkoff Theorem $p-\lambda q$ is strictly hyperbolic and of 
degree at least $n-1$. Clearly, $T(p-\lambda q)=0$,  
which proves the Claim.  

Suppose now that $T$ maps all strictly stable polynomials of degree $n$ into
the set $\HH_1(\CC)\cup\{0\}$ (the case when $T$ maps all such polynomials
into $\HH_1^-(\CC) \cup\{0\}$ is treated similarly). Let $f$ be a strictly 
stable polynomial with $\deg f<n$ and set 
$f_{\eps}(z)=(1-\eps iz)^{n-\deg f}f(z)$. Then $f_{\eps}$ is strictly stable 
whenever $\eps>0$ and $T(f_{\eps})\in
\HH_1(\CC)\cup\{0\}$ since $\deg f_{\eps}=n$. Letting $\eps\to 0$ we get
$T(f)\in \HH_1(\CC)\cup\{0\}$ and since strictly stable polynomials are 
dense in $\HH_1(\CC)$ it follows that $T$ is stability preserving.
\end{proof}

As a final tool we will need the Grace-Walsh-Szeg\"o coincidence theorem 
\cite{grace,szego,walsh}.  
Recall that a multivariate polynomial is {\em multi-affine} if it has degree 
at most 
one in each variable. 

\begin{theorem}[Grace-Walsh-Szeg\"o coincidence theorem]
Let $f \in \CC[z_1,\ldots, z_n]$ be symmetric and multi-affine and let $C$ be 
a circular domain 
containing the points $\zeta_1, \ldots, \zeta_n$. Suppose that either the 
total degree of $f$ equals $n$ or that $C$ is convex (or both). 
Then there exists at least one point $\zeta \in C$ such that 
$$
f(\zeta_1, \ldots, \zeta_n)=f(\zeta, \ldots, \zeta).
$$
\end{theorem}

From the Grace-Walsh-Szeg\"o coincidence theorem we immediately deduce:

\begin{corollary}\label{polarizeit}
Let $f \in \CC[z_1,\ldots, z_n]$ be of degree at most  $d$ in $z_1$ and 
consider the expansion of $f$ in powers of $z_1$, 
$$
f(z_1,\ldots,z_n)= \sum_{k=0}^dQ_k(z_2,\ldots, z_n)z_1^k,\quad 
Q_k\in\bC[z_2,\ldots,z_n],\, 0\le k\le d. 
$$
Then $f$ is stable if and only if the polynomial 
$$
\sum_{k=0}^d Q_k(z_2,\ldots,z_n)\frac {e_k(x_1,\ldots, x_d)}{\binom d k}
$$
is stable in the variables $z_2, \ldots, z_n, x_1, \ldots, x_d$, where 
$e_k(x_1, \ldots, x_d)$, 
$0 \leq k \leq d$, are the elementary symmetric functions in the variables 
$x_1, \ldots, x_d$ given 
by $e_0=1$ and 
$e_k = \sum_{1 \leq j_1 < j_2 < \cdots < j_k \leq d} x_{j_1}\cdots x_{j_k}$ 
for 
$1 \leq k \leq d$. 
\end{corollary}

\begin{lemma}\label{extension}
Let $T: \CC_n[z_1] \rightarrow \CC[z_1]$ be a linear operator such that 
$T[(z_1+w)^n] \in \HH_2(\CC)$. If $f \in \HH_m(\CC)$ is of degree at most $n$ 
in $z_1$ then 
$T(f) \in \HH_m(\CC)\cup\{0\}$, where $T$ is extended to a linear operator on
$\bC[z_1,\ldots,z_m]$ by setting 
$
T(z_1^{\alpha_1}\cdots  z_m^{\alpha_m}) 
= T(z_1^{\alpha_1})z_2^{\alpha_2} \cdots z_m^{\alpha_m}
$ for all $\alpha \in \NN^m$ with $\alpha_1 \leq n$ (compare with 
Notation~\ref{not1}).
\end{lemma}

\begin{proof}
Let $f \in \HH_m(\CC)$ be of degree $n$ in $z_1$. For $\epsilon >0$ set 
$$
f_\epsilon(z_1, \ldots, z_m) = 
f(z_1+\epsilon i, z_2,\ldots,z_m).
$$ 
Fixing $\zeta_2, \ldots, \zeta_m$ in the open upper half-plane  
we may write $f_\epsilon(z_1, \zeta_2, \ldots, \zeta_m)$ as
$$
f_\epsilon(z_1, \zeta_2, \ldots, \zeta_m)
= C(z_1-\xi_1)(z_1-\xi_2)\cdots (z_1-\xi_n) = 
C \sum_{k=0}^n(-1)^k e_k(\xi_1,\ldots, \xi_n)z_1^{n-k}, 
$$
where $C \neq 0$ and $\Im(\xi_j)<0$, $1\le j\le n$. 
 Note that $z_1\mapsto f_\epsilon(z_1, \zeta_2, \ldots, \zeta_m)$ is indeed a 
polynomial of degree $n$ in $z_1$. This is because 
$$f_\epsilon(z_1, z_2, \ldots, z_m)= z_1^nQ_n(z_2, \ldots, z_m)
+\mbox{terms of lower degree in $z_1$}$$
and $Q_n(z_2, \ldots, z_n) := 
\lim_{r \rightarrow \infty}r^{-n}f(r,z_2,\ldots,z_m)$ is stable by 
Hurwitz' theorem hence $Q_n(\zeta_2,\ldots,\zeta_n)\neq 0$.
Let us now write $T[(z_1+w)^n]$ as 
$$
T[(z_1+w)^n] = \sum_{k=0}^n \binom n k T\left(z_1^{n-k}\right)w^k \in
 \HH_2(\CC).
$$
By Corollary \ref{polarizeit} we know that 
$$
 \sum_{k=0}^n T\left(z_1^{n-k}\right)e_k(w_1,\ldots, w_n) \in \HH_{n+1}(\CC).
$$ 
But then 
$$
T(f_\epsilon)(z_1, \zeta_2, \ldots, \zeta_m)=  C\sum_{k=0}^n 
T\left(z_1^{n-k}\right)e_k(-\xi_1,\ldots, -\xi_n) \in \HH_{1}(\CC), 
$$ 
which gives $T(f_\epsilon) \in \HH_m(\CC)$. Letting $\epsilon \rightarrow 0$ 
we have $T(f) \in \HH_m(\CC)\cup\{0\}$. 
If $f \in \HH_m(\CC)$ is of degree less than $n$ in $z_1$ we may consider 
$f^\epsilon=(1-\epsilon i z_1)^{n-\deg f} f \in \HH_m(\CC)$.  Then by the 
above  
one has $T(f^\epsilon) \in  \HH_n(\CC)\cup\{0\}$ for all $\epsilon>0$. The 
lemma now follows from Hurwitz' theorem by letting $\epsilon \rightarrow 0$. 
\end{proof}

We are now ready to prove Theorem \ref{finitestab}. 
\begin{proof}[Proof of Theorem \ref{finitestab}]
If $T[(z+w)^n] \in \HH_2(\CC)$ then by applying Lemma \ref{extension} with 
$m=1$ it follows that $T$ is stability 
preserving.

Suppose now that $T$ preserves stability. Assume first that there exists 
$w_0\in\bC$ with $\Im(w_0)>0$ such that $(z+w_0)^n$ is in the kernel of $T$. 
Since $(z+w_0)^n$ is strictly stable it follows from  Lemma \ref{degenerate} 
and Lemma \ref{spaces}  that 
$\dim_{\bC} T(\CC_n[z]) \leq 1$. Hence $T$ is given by $T(f) = \alpha(f)P$, 
where $P$ is a fixed stable 
polynomial and $\alpha : \CC_n[z]  \rightarrow \CC$ is a linear functional. 
Otherwise we may assume that  $T[(z+w_0)^n] \in \HH_1(\CC)$ for all 
$w_0 \in \CC$ with 
$\Im(w_0) >0$ and conclude that $T[(z+w)^n] \in \HH_2(\CC)$. 
\end{proof}

The proof of Theorem \ref{finitehyp} now follows easily. 

\begin{proof}[Proof of Theorem \ref{finitehyp}]
Recall Notation~\ref{not4} and note that by Lemma \ref{trichotomy} we may 
assume that $\dim_{\bR} T(\RR_n[z]) > 2$,  so that 
$T$ is either stability reversing or stability preserving. By 
Theorem \ref{finitestab} $T$ is 
stability reversing or stability preserving if and only if 
$T[(z+w)^n] \in \HH_2(\CC)$ or 
$T[(z-w)^n] \in \HH_2(\CC)$, respectively, which proves the theorem. 
\end{proof}

\begin{proof}[Proof of Theorem \ref{finitehypC}]
By Theorem \ref{finitehyp} it suffices to prove that if a linear operator  
$T : \CC_n[z] \rightarrow \CC[z]$ 
satisfies $T : \pi_n(\RR) \rightarrow \pi(\RR)$ then either 
\begin{itemize}
\item[(a)] There exists $\theta \in \RR$ such that $T=e^{i\theta}\tilde{T}$, 
where $\tilde{T}: \RR_n[z] \rightarrow \RR[z]$ is 
a hyperbolicity preserver when restricted to $\RR_n[z]$, or 
\item[(b)] $T$ is given by $T(f) = \alpha(f)P$, where 
$\alpha : \CC_n[z] \rightarrow \CC$ is a linear functional and $P$ is a 
hyperbolic polynomial. 
\end{itemize}
We prove this using induction on $n \in \NN$.  
If $n=0$ there is nothing to prove so we may assume that $n$ is a positive 
integer. Note that $T$ restricts to 
a linear operator $T' : \pi_{n-1}(\RR) \rightarrow \pi(\RR)$ and that by 
induction $T'$ must be of the form (a) or (b) above. Suppose that 
$T'=e^{i\theta}\tilde{T}'$, where $\tilde{T}'$ is a hyperbolicity preserver 
up to degree $n-1$. If $T(z^n)=e^{i\theta}f_n$, where $f_n \in \RR[z]$ 
(actually, $f_n\in\HH_1(\bR)\cup\{0\}$), then $T$ is of the form (a). Hence 
we may assume that $T(z^n)= e^{i\gamma}f_n$, where 
$0 \leq \gamma < 2\pi$, $\gamma- \theta$ is not an integer multiple of 
$\pi$ and 
$f_n$ is a hyperbolic polynomial. Suppose that there is an integer $k < n$ 
such that $T(z^k)$ is not a 
constant multiple of $f_n$. Let $M$ be the largest such $k$ and set 
$R(z)=e^{-i\theta}T(z^M)$. Then 
$$
e^{-i\theta}T\!\left[z^M(1+z)^{n-M}\right] 
= R(z)+ \left(r+e^{i(\gamma-\theta)}\right)f_n(z)
$$
for some  $r \in \RR$. But this polynomial is supposed to be a complex 
constant multiple of 
a hyperbolic polynomial, which can only happen if $R$ and thus 
$T\left(z^M\right)$ is a constant multiple of $f_n$.
This contradiction means that $T$ must be as in (b) above with $P=f_n$. 

Assume now that $T'$ is as in (b). 
If $T(z^n)$ is a constant multiple of $P$ there is nothing to prove, so we 
may assume that $T(z^n)=e^{i\theta}f_n$, where $0 \leq \theta < 2\pi$ and 
$f_n$ is a hyperbolic polynomial which is not a constant multiple of $P$. 
Suppose that there is an integer $k <n$ such 
that $\alpha(z^k)e^{-i\theta} \notin \RR$ and let $M$ be the largest such 
integer. Then 
$$
e^{-i\theta}T\!\left[z^M(1+z)^{n-M}\right]
= \left(\alpha\!\left(z^M\right)e^{-i\theta}+r\right)P(z) + f_n(z)
$$
for some $r \in \RR$. However, the latter polynomial should be a complex 
constant multiple of 
a hyperbolic polynomial and this can happen only if  
$\alpha(z^M)e^{-i\theta} \in \RR$, which 
contradicts the above assumption. This means that $T$ must be as in (a). 
\end{proof}

\begin{proof}[Proof of Corollary~\ref{cor-alghyp}]
Note first that if $T:\bR[z]\to \bR[z]$ is as in (a), (b) or (c) of 
Corollary \ref{cor-alghyp} then by Theorem~\ref{finitehyp} we have that $T$ 
preserves hyperbolicity up to any degree $n\in\bN$. Conversely, if 
$T:\bR[z]\to \bR[z]$ preserves hyperbolicity then for any $n\in\bN$ the 
restriction $T:\bR_n[z]\to \bR[z]$
preservers hyperbolicity (up to degree $n$). The case when 
$\dim_{\RR} T(\RR[z])\le 2$ is clear. Suppose now that 
$\dim_{\RR} T(\RR[z])>2$ and that 
$T[(z+w)^n] \in \HH_2(\RR)$ for some $n\in\bN$. Then by 
Lemma \ref{extension} we have that 
$T[(z+w)^m] \in \HH_2(\CC)\cup\{0\}$ for all $m \leq n$ and since the latter 
polynomial has real coefficients we get $T[(z+w)^m] \in \HH_2(\RR)\cup\{0\}$ 
for all $m \leq n$. Similarly, if $T[(z-w)^n] \in \HH_2(\RR)$ then 
$T[(z-w)^m] \in \HH_2(\RR)\cup\{0\}$ for all $m \leq n$. It follows that if 
$\dim_{\RR} T(\RR[z])>2$ then 
either  $T[(z+w)^n] \in \HH_2(\RR)\cup\{0\}$ for all $n \in \NN$ or 
$T[(z-w)^n] \in \HH_2(\RR)\cup\{0\}$ for all $n \in \NN$.
\end{proof}

By appropriately modifying the arguments in the proof of 
Corollary~\ref{cor-alghyp} one can easily check that 
Corollary~\ref{cor-algstab} follows readily from Theorem~\ref{finitestab}.

\subsection{The Transcendental Characterization}\label{ss32}

We will need the following lemma due to 
Sz\'asz, see \cite[Lemma 3]{szasz}. 
\begin{lemma}[Sz\'asz]\label{szasz} 
Let $m,n\in\bN$ with $m\le n$ and $f(z)=\sum_{k=m}^{n}c_kz^k\in\bC[z]$. If
$f(z)\in\HH_1(\CC)$ and $c_mc_n \neq 0$ then for any $r\ge 0$ one has  
$$
|f(z)| \leq |c_m|r^m\exp \left(r\frac{|c_{m+1}|}{|c_m|}
+3r^2\frac{|c_{m+1}|^2}{|c_m|^2}+3r^2\frac{|c_{m+2}|}{|c_m|}\right) 
$$
whenever $|z| \leq r$. 
\end{lemma} 

For $k,n \in \NN$ let $(n)_k = k!\binom n k = n(n-1)\cdots (n-k+1)$ if 
$k \leq n$ and $(n)_k =0$ if $n<k$ denote as usual the Pochhammer symbol.  

\begin{theorem}\label{genpol}
Let $F(z,w)=\sum_{k=0}^{\infty}P_k(z)w^k$ be a formal power series in $w$ 
with polynomial coefficients. Then $F(z,w) \in \HC_2(\CC)$ if and only if 
$\sum_{k=0}^n(n)_ k P_k(z) w^k \in \HH_2(\CC)\cup\{0\}$ for 
all $n \in \NN$. 
\end{theorem}

\begin{proof}
Suppose that  $F(z,w)=\sum_{k=0}^\infty P_{k}(z)w^k\in \HC_2(\CC)$ has 
polynomial coefficients. Given $n \in \NN$, the sequence $\{(n)_k\}_{k=0}^n$ 
is a 
multiplier sequence, and since it is non-negative it is a stability preserver 
by Theorem~\ref{trichotomy}.  By   
Corollary~\ref{cor-algstab} and Lemma~\ref{extension} we have that this 
multiplier extends 
to a map $\Lambda : \HH_2(\CC) \rightarrow \HH_2(\CC)\cup\{0\}$. Now, if  
$\tilde{F}_m(z,w) = \sum_{k=0}^{N_m}P_{m,k}(z)w^k$ is a sequence of
 polynomials in $\HH_2(\CC)$ converging to $F(z,w)$, uniformly on compacts, 
then 
$P_{m,k}(z) \rightarrow P_k(z)$ as $m \rightarrow \infty$  uniformly on 
compacts for fixed $k \in \NN$. But then  
$\Lambda[\tilde{F}_m(z,w)] \rightarrow \Lambda[F(z,w)]$ uniformly on 
compacts, which gives $\sum_{k=0}^n(n)_ k P_k(z) w^k \in \HH_2(\CC)\cup\{0\}$.

Conversely, suppose that 
$\sum_{k=0}^n(n)_ k P_k(z) w^k \in \HH_2(\CC)\cup\{0\}$ for all $n \in \NN$ 
and let 
$F_n(z,w) = \sum_{k=0}^n(n)_ k n^{-k}P_k(z) w^k$. 

\medskip

{\em Claim.} Given $r>0$ there is a constant $C_r$ such that 
$$
|F_n(z,w)| \leq C_r \mbox{ for } |z| \leq r, |w| \leq r 
\mbox{ and all } n \in \NN. 
$$

This claim proves the theorem since $\{ F_n(z,w) \}_{n=0}^\infty$ is then a 
normal family whose 
convergent subsequences converge to $F(z,w)$ (by the fact that 
$n^{-k}(n)_k\rightarrow 1$ for all $k \in \NN$ as $n \rightarrow \infty$). 

We first prove the Claim in the special case 
when $P_k(z) \in \RR[z]$ for all $k\in \NN$ and $P_K(z)$ is a non-zero 
constant, where $K$ is the first index 
for which $P_K(z) \neq 0$. 

\medskip

{\em Proof of the Claim in the special case.} Let $|P_K(z)|=A$, 
$B_r =\max\{ |P_{K+1}(z)| : |z|\leq r\}$ and 
$D_r =\max\{ |P_{K+2}(z)| : |z|\leq r\}$. Then, 
if we fix $\zeta \in \CC$ with $\Im(\zeta)>0$, we have that 
$F_n(\zeta, w) \in \HH_1(\CC) \cup\{0\}$ and by Lemma \ref{szasz} 
\begin{equation}\label{uppskatta}
|F_n(\zeta,w)| \leq  
Ar^K\exp\!\left(r\frac{B_r}{A}+3r^2\frac{B_r^2}{A^2}+3r^2\frac{D_r}{A}\right)  
\end{equation}
whenever $\Im(\zeta)>0$, $|\zeta| \leq r$ and $|w| \leq r$. If 
$\zeta \in \CC$ is fixed with  $\Im(\zeta)<0$ then $F_n(\zeta,-w) 
\in \HH_1(\CC)$ (since 
$F_n(z,w)$ has real coefficients and 
$\overline{F_n(z,w)}=F_n(\overline{z}, \overline{w})$).  
By Lemma \ref{szasz} this means that \eqref{uppskatta} holds also 
for $\Im(\zeta)<0$ and by continuity also for $\zeta \in \RR$, which proves 
the Claim. 

Next we assume that $\deg(P_K(z)) =d \geq 1$. An application of 
Theorem \ref{finitestab} verifies that $T = d/dz$ preserves stability, and 
by Lemma \ref{extension} $T= \partial/\partial z$ preserves stability in two 
variables. Hence 
$\frac{\partial F_n(z,w)}{\partial z} \in \HH_2(\RR)\cup\{0\}$ if 
$F_n(z,w) \in \HH_2(\RR)\cup\{0\}$. To deal with this case it is therefore 
enough to prove that if 
$\left|\frac{\partial F_n(z,w)}{\partial z}\right| \leq C_r$ for 
$|z| \leq r, |w| \leq r$ and all $n \in \NN$ then there 
is a constant $D_r$ such that $|F_n(z,w)| \leq D_r$ for 
$|z| \leq r, |w| \leq r$ and all $n \in \NN$. 
Clearly, $F_n(0,w) \in \HH_1(\RR) \cup \{0\}$ for all $n\in \NN$. 
Moreover, if $m$ is the first index 
such that $P_m(0) \neq 0$ then for $n \geq m$ we have 
$$
F_n(0,w)=
(n)_mn^{-m}P_m(0)w^m +(n)_{m+1}n^{-m-1}P_{m+1}(0)w^{m+1}+\ldots\in \HH_1(\RR),
$$
which by Lemma \ref{szasz} gives  
$$
|F_n(0,w)| \leq  |P_m(0)|r^m\exp \left(r\frac{|P_{m+1}(0)|}{|P_{m}(0)|}
+3r^2\frac{|P_{m+1}(0)|^2}{|P_{m}(0)|^2}
+3r^2\frac{|P_{m+2}(0)|}{|P_{m}(0)|}\right) =: E_r
$$
for $|z| \leq r$. Here we have used that $(n)_kn^{-k} \leq 1$ and 
$(n)_{k+1}n^{-k-1}/(n)_kn^{-k} \leq 1$  for $0 \leq k \leq n$. 
Suppose that 
$$
\left|\frac{\partial F_n(z,w)}{\partial z}\right| \leq C_r
$$
for $n \in \NN$ and $|z|\leq r$, $|w|\leq r$. Then, 
$$
F_n(z,w) = F_n(0,w)+ z\int_0^1 \frac{\partial F_n}{\partial z}(zt,w) dt,
$$
so 
$$
|F_n(z,w)| \leq E_r + rC_r \ \  \mbox{ for } |z| \leq r, 
|w| \leq r \mbox{ and all } n \in \NN. 
$$

Next we prove the above Claim in the general case. For this we need a 
property of multivariate stable polynomials that was established in 
\cite[Corollary 1]{BBS}: 

\medskip

{\em Fact.} If $h =g+if \in \HH_n(\CC)$ then $f,g \in \HH_n(\RR)$. 

\medskip

Now this means that if we let $P_k(z)= R_k(z)+i I_k(z)$ then we may write 
$$F_n(z,w) = F^{\Re}_n(z,w) + iF^{\Im}_n(z,w)$$ 
with $F^{\Re}_n(z,w), F^{\Im}_n(z,w) \in \HH_2(\RR)\cup\{0\}$, where 
$F^{\Re}_n(z,w), F^{\Im}_n(z,w)$ are given by 
$$F^{\Re}_n(z,w) = \sum_{k=0}^n(n)_ k n^{-k}R_k(z) w^k,\quad 
F^{\Im}_n(z,w) = \sum_{k=0}^n(n)_ k n^{-k}I_k(z) w^k.$$ 
By the above there are constants $A_r$ and $B_r$ such that 
$$|F^{\Re}_n(z,w)| \leq A_r\text{ and }|F^{\Im}_n(z,w)| \leq B_r
\text{ for }|z| \leq r, |w| \leq r\text{ and all }n \in \NN.$$ 
Hence  $|F_n(z,w)| \leq \sqrt{ A_r^2+B_r^2}$ 
for $|z| \leq r, |w| \leq r$ and all $n \in \NN$, which proves the Claim in 
the general case.
\end{proof}

We can now prove the transcendental characterizations of hyperbolicity and 
stability preservers, respectively. 

\begin{proof}[Proof of Theorems \ref{tranR} and \ref{tranC}]
We only prove Theorem~\ref{tranC} since the proof of Theorem~\ref{tranR} is 
almost identical. 
Theorem~\ref{tranC} follows quite easily from 
Theorem~\ref{genpol} and Corollary~\ref{cor-algstab}. Indeed, note that 
$T[(1-zw)^n] \in \HH_2(\CC)$ if and only if 
$T[(z+w)^n] \in \HH_2(\CC)$. Since
$$
T[(1-zw)^n] = \sum_{k=0}^n (n)_k \frac{(-1)^kT(z^k)}{k!} w^k 
$$
for all $n\in\bN$ we deduce the desired conclusion by comparing the above 
expression with the modified symbol $G_T(z,w)$ introduced in 
Notation~\ref{not2}.
\end{proof}

Finally, we show how P\'olya-Schur's algebraic and transcendental 
characterizations of multiplier sequences (Theorem~\ref{ps}) 
follow from our results. 

\begin{proof}[Proof of Theorem \ref{ps}]
We may assume that  $\lambda(0)=1$. Statements (i)--(iv) are all true 
if $\dim_{\bR}T(\RR[z]) \leq 2$ since 
then $\lambda(j)=0$ for all $j \geq 2$. Hence we may assume that 
$\dim_{\bR}T(\RR[z]) >2$. 
Corollary \ref{cor-alghyp} implies that either 
$$
T[(1-zw)^n]=\sum_{k=0}^n(-1)^k\binom n k \lambda(k) (zw)^k  \in \HH_2(\RR)
$$
for all $n \in \NN$ or  
$$
T[(1+zw)^n]  = \sum_{k=0}^n \binom n k \lambda(k) (zw)^k \in \HH_2(\RR)
$$
for all $n \in \NN$. We claim that if  $f(z) \in \RR[z]$ then
$
f(zw) \in \HH_2(\RR) 
$
if and only if all the zeros of $f$ are real and non-negative. Suppose first 
that $f(zw)$ is real stable. Letting $z=w=t$ we see that 
$f(t^2)$ is hyperbolic, which can only happen if all the zeros of $f$ are 
non-negative (since $a+t^2$ is hyperbolic if and only if $a \leq 0$). On the 
other hand, if $f$ has only real non-positive zeros then $f(zw)$ factors as   
$f(zw)=C\prod_{j=0}^n(zw+\alpha_j),$
where $\alpha_j \in \RR$, $1 \leq j \leq n$. 
Now $zw+\alpha_j \in \HH_2(\RR)$ if and only if $\alpha_j \leq 0$. Hence (i) 
$\Leftrightarrow$ (iv) and by 
Theorem \ref{tranR} we also have (ii)  $\Leftrightarrow$ (iv). The 
equivalence of (ii) and (iii) is Laguerre's 
classical result. 
\end{proof}

\subsection{Closed Circular Domains and Their Boundaries}\label{ss33}

Recall Definitions~\ref{def1} and~\ref{def-circ}, Notation~\ref{not-sets} and 
the linear operators
$\phi_n$ introduced in Notation~\ref{not-fi}. In particular, an $H$-stable 
polynomial in the sense of Definition~\ref{def-circ} is precisely a stable
polynomial in the sense of Definition~\ref{def1}.

\begin{lemma}\label{save}
Let $T : \CC_n[z] \rightarrow \CC_m[z]$ be a linear operator and suppose $m$ 
is minimal, i.e., $m= \max\{\deg T(f): f \in \CC_n[z]\}$. Let further 
$C=\Phi^{-1}(H)$ be an open circular domain, where $\Phi$ is a M\"obius 
transformation as in \eqref{mobius}, and let 
$S : \CC_n[z] \rightarrow \CC_m[z]$  be the linear operator defined by 
$S= \phi^{-1}_m T \phi_n$. 
The following are equivalent: 
\begin{itemize}
\item[(i)] $T(f)$ is $C$-stable or zero whenever $f$ is of degree $n$ and 
$C$-stable, 
\item[(ii)] $S(f)$ is $H$-stable or zero whenever $f$ is of degree $n$ and 
$H$-stable,  
\item[(iii)] $S(f)$ is $H$-stable or zero whenever $f$ is of degree at most 
$n$ and $H$-stable.  
\end{itemize}
The following are also equivalent: 
\begin{itemize}
\item[(iv)] $T(f)$ is $\partial C$-stable or zero whenever $f$ is of degree 
$n$ and $\partial C$-stable, 
\item[(v)] $S(f)$ is $\RR$-stable or zero whenever $f$ is of degree $n$ and 
$\RR$-stable,  
\item[(vi)] $S(f)$ is $\RR$-stable or zero whenever $f$ is of degree at most 
$n$ and $\RR$-stable. 
\end{itemize}
\end{lemma}

\begin{proof}
Note first that the equivalences (ii)$\Leftrightarrow$ (iii) and 
(v)$\Leftrightarrow$ (vi) are simple consequences of the density argument 
used in \S \ref{ss22}. Let us now show that (i) $\Leftrightarrow$ (ii). 
This  is obvious if $C$ is an open half-plane, i.e., if $c=0$ 
(cf.~\eqref{assum}). Therefore we assume that $c \neq 0$ and 
that the boundary of $C$ is a circle. If $C$ is an open disk then $a/c$ is in 
the open lower half-plane, so that $-cz+a$ is $H$-stable. Moreover, 
$-d/c\in\partial C$ and thus $cz+d$ is $C$-stable. Hence $f$ is $H$-stable 
if and only if $\phi_n(f)$ is $C$-stable so the assertion follows in this 
case as well. 

It remains to prove that (i) $\Leftrightarrow$ (ii) in the case when $C$ is 
the open complement  of a (closed) disk, which we proceed to do. 

\medskip

(i) $\Rightarrow$ (ii). Clearly, we may assume that $T$ is not the trivial 
(identically zero) operator. 
Let $p(z)=\sum_{k=0}^na_kz^k$ be an $H$-stable polynomial of degree $n$ and 
suppose first that $\sum_{k=0}^na_ka^kc^{n-k}\neq 0$. Then $\phi_n(p)$ is a 
$C$-stable polynomial of degree   
$n$ so that by assumption $T(\phi_n(p))$ is $C$-stable or zero. If 
$S(p) \neq 0$ it follows that 
that we can express $S(p)$ uniquely as 
$S(p) = (-cz+a)^{r(p)}S_0(p)$,  
where $S_0(p)$ is $H$-stable and $r(p)$ is a non-negative integer. By a 
continuity argument and an application of Hurwitz' theorem we have that a 
factorization as above holds for any $H$-stable polynomial of degree $n$.
 Since the set of $H$-stable polynomials  
of degree $n$ is dense in $\pi_n(H')\cup \{0\}$ -- that is, the set of 
$H$-stable polynomials of degree at most $n$ union the zero 
polynomial -- we deduce that such a 
factorization holds for the image under $S$ of any $H$-stable polynomial $p$ 
of degree at most $n$, namely 
$$S(p) = (-cz+a)^{r(p)}S_0(p),$$   
where $S_0(p)$ is $H$-stable or zero and $r(p)$ is a non-negative integer.

Fix a basis $\{p_j(z)\}_{j=0}^n$ of  $\CC_n[z]$ consisting of 
{\em strictly $H$-stable} (that is, $\overline{H}$-stable) polynomials of 
degree $n$. We distinguish two cases:

Suppose first that $S(f) \neq 0$ for all strictly $H$-stable polynomials $f$ 
of degree $n$. Since the topological space of strictly stable polynomials of 
degree $n$ is (path-) connected we  have by Hurwitz' theorem that $r(f)$ is 
constant on the set  of strictly $H$-stable polynomials of degree $n$. 
Thus, by the minimality assumption on $m$ and the fact that $\phi_n$ is 
invertible we must 
have $\deg(T(\phi_n(p_k)))=m$ for some $k$. It follows that $r(p_k)=0$ hence
$S(f)=S_0(f)$ for any strictly $H$-stable polynomial $f$ of degree $n$. 
Using again a standard density argument and Hurwitz' theorem we deduce that 
$S$ preserves ($H$-)stability up to degree $n$. 

Suppose now that  $S(f) = 0$ for some strictly $H$-stable polynomial $f$ of 
degree $n$ and let $g \in \CC_n[z]$. Clearly, $f + \epsilon g$ is strictly 
$H$-stable for all $\epsilon>0$ small enough. By the above we have that 
$$
S(g)=\epsilon^{-1}S(f+\epsilon g)= (-cz+a)^{r(g)}S_0(g),
$$
where $S_0(g)$ is $H$-stable or zero. It follows that $V:=S(\CC_n[z])$ is a 
$\CC$-linear space such that every non-zero polynomial in $V$ is a 
$(-cz+a)^r$-multiple of an $H$-stable polynomial. We know that $r(p_k)=0$ 
for the strictly $H$-stable polynomial $p_k$ above. Assume that 
$h \in \CC_n[z]$ is such that $r(h)\neq 0$. Since $r(p_k)=0$ and 
$$
S(h) + \delta S(p_k)=S(h+\delta p_k)
=(-cz+a)^{r(h+\delta p_k)}S_0(h+\delta p_k)\in V,
$$ 
where $S_0(h+\delta p_k)$ is $H$-stable or zero, we conclude that 
$S(h) + \delta S(p_k)$ is $H$-stable or zero for all $\delta \neq 0$. 
Letting $\delta \rightarrow 0$ we have by Hurwitz' theorem that either 
$S(h)=0$ or $S(h)$ is $H$-stable. However, this contradicts the fact that 
$a/c\in H$ and $S(h)(a/c)=0$, which follows from the assumption that 
$r(h)\neq 0$. Hence all non-zero polynomials in $V$ are $H$-stable and thus 
we deduce that $S$ 
preserves ($H$-)stability up to degree $n$ in this case as well. 

\medskip

(ii) $\Rightarrow$ (i). Since the set of $H$-stable polynomials of degree $n$ 
is dense in the set of $H$-stable polynomials of degree at most $n$ it 
follows that $S$ preserves $H$-stability on the latter set 
(i.e., $S$ preserves $H$-stability up to degree $n$). Let $f$ be a 
$C$-stable polynomial of degree $n$. Then $\phi_n^{-1}(f)$ is $H$-stable. 
Hence so is $S(\phi_n^{-1}(f))$ and thus 
$T(f)=\phi_m\big(S(\phi_n^{-1}(f))\big)$ is a $C$-stable polynomial (note 
that $-d/c\in\partial C$). 

\medskip

The equivalence (iv) $\Leftrightarrow$ (v) follows just as above by 
replacing ``strictly $H$-stable'' with ``strictly hyperbolic'', that is, 
real- and simple-rooted. 
\end{proof}

\begin{notation}\label{not-fifi}
Given a polynomial $f \in \CC[z,w]$ 
of degree at most $m$ in $z$ and at most $n$ in $w$ and a M\"obius 
transformation $\Phi$ as in \eqref{mobius} we let 
$$\phi_{m,z}(f)(z,w)= (cz+d)^mf(\Phi(z),w),\quad 
\phi_{n,w}(f)(z,w)= (cw+d)^nf(z,\Phi(w)).$$
\end{notation}

\begin{lemma}\label{fifi} 
Let $f(z,w) \in \CC[z,w]$ be of degree at most  $m$ in $z$ and at most $n$ in 
$w$ and let $\Phi : C \rightarrow H$ be a M\"obius transformation  as 
in \eqref{mobius}. If either
\begin{itemize}
\item[(a)] $C$ is not the exterior of a disk, or 
\item[(b)] $C$ is the exterior of a disk and 
\begin{itemize} 
\item[(b1)] the degree in $z$ of $\phi_{m,z}(f)(z,w)$ is $m$, and 
\item[(b2)] the degree in $w$ of $\phi_{n,w} \phi_{m,z}(f)(z,w)$ is $n$, 
\end{itemize} 
\end{itemize}
then $f$ is $H$-stable if and only if $\phi_{n,w}\phi_{m,z}(f)$ is $C$-stable.
\end{lemma}

\begin{proof}
The equivalence is clear when $C$ is a disk or a half-plane since in these 
cases $-cz+a$ is $H$-stable and $cz+d$ is $C$-stable (cf.~\eqref{assum}). 
Hence we may assume that $C$ is the exterior of a disk. Note that if $g(z)$ 
is a polynomial of degree $k$ then $a/c$ is not a zero of 
$\phi_k^{-1}(g)(z)$. In particular, if  $g(z)$ is a $C$-stable polynomial of 
degree $k$ then $\phi_k^{-1}(g)(z)$ is an $H$-stable polynomial of 
degree $k$.  Now since $-d/c \in \partial C$ we have that 
$cz+d$ is $C$-stable and therefore $\phi_{n,w} \phi_{m,z}(f)$ is $C$-stable 
if $f$ is $H$-stable, which proves one of the implications. 

Conversely, suppose that 
$G(z,w):=\phi_{n,w} \phi_{m,z}(f)(z,w)= \sum_{k=0}^nQ_k(z)w^k$ 
is $C$-stable. Then so is 
$\lambda^{-n}G(z,\zeta +\lambda(w-\zeta))$ whenever $\lambda \geq 1$, where 
$\zeta$ is the center of 
$C'$. Letting $\lambda \rightarrow \infty$ we see by Hurwitz' theorem that 
$(w-\zeta)^nQ_n(z)$ is $C$-stable. Hence, so is $Q_n(z)$. For every 
$z_0 \in C$ the degree of $G(z_0,w)$ is $n$, from which it follows that 
$\phi^{-1}_{n,w}(G)(z,w)=\phi_{m,z}(f)(z,w) \neq 0$ whenever 
$z \in C$, $w \in H$. Similarly, 
if $\phi_{m,z}(f)(z,w)= \sum_{k=0}^m P_k(w)z^k$ then $P_m(w)$ is $H$-stable 
so that 
$\phi_{m,z}(f)(z,w_0)$ has degree $m$ for every $w_0 \in H$. We conclude that 
$f(z,w)=\phi^{-1}_{m,z}\phi^{-1}_{n,w}(G)(z,w)$ is $H$-stable, as claimed.
\end{proof}

The following lemma is a simple consequence of Hurwitz' theorem since 
boundedness prevents zeros from escaping to infinity. 

\begin{lemma}\label{m}
Let $\Omega_1$ be a path-connected subset of $\CC$ and let $\Omega_2$ be a 
bounded subset of $\CC$. If $T : \CC_n[z] \rightarrow \CC[z]$ is a linear 
operator 
such that 
$T : \pi_n(\Omega_1)\setminus \pi_{n-1}(\Omega_1) \rightarrow \pi(\Omega_2)$ 
then 
all polynomials in the image of  
$\pi_n(\Omega_1)\setminus \pi_{n-1}(\Omega_1)$ have the same degree. 
\end{lemma}

Note that in the hypothesis of Lemma~\ref{m} we do not allow the identically 
zero polynomial to be in the image of 
$\pi_n(\Omega_1)\setminus \pi_{n-1}(\Omega_1)$. 

\begin{proof}[Proof of Theorem~\ref{circularregion}]
Suppose that $T : \CC_n[z] \rightarrow \CC_m[z]$, where as before $m$ is  
minimal in the sense that $m = \max\{\deg T(f): f \in \CC_n[z]\}$. By 
combining Lemma~\ref{save} with Hurwitz' theorem we see that 
$T : \pi_n(C')\setminus \pi_{n-1}(C') \rightarrow \pi(C')\cup\{0\}$ if 
and only if 
$\phi^{-1}_mT\phi_n : \pi_n(H') \rightarrow \pi(H')\cup\{0\}$.  
The case  $\dim_{\bC} T(\CC_n[z]) \leq1$ is clear. 
If $\dim_{\bC} T(\CC_n[z]) >1$ then by Theorem \ref{finitestab} we 
have $\phi_m^{-1}T\phi_n : \pi_n(H') \rightarrow \pi(H')\cup\{0\}$ if and only 
if $f(z,w):=\phi^{-1}_{m}T\phi_n\left[(z+w)^n\right]$ is $H$-stable. 
Now the polynomial in (b) of Theorem~\ref{circularregion} is precisely 
$\phi_{n,w} \phi_{m,z}(f)$. Moreover, by Lemma~\ref{fifi} (a) we may assume 
that $C$ is the exterior of a disk. Therefore, in order to complete the 
proof it only remains to show that conditions (b1) and (b2) of 
Lemma~\ref{fifi} are satisfied. 

If there were two different  $\overline{C}$-stable polynomials of degree 
$n$ that were mapped by $T$ on polynomials of different degrees then by 
Lemma~\ref{m} there would be a  $\overline{C}$-stable polynomial $g$ of 
degree $n$ in the kernel of $T$. However, since $\phi^{-1}_n(g)$ is 
strictly $H$-stable it would then 
follow from Lemmas \ref{degenerate} and \ref{spaces} 
that $\phi_m^{-1}T\phi_n$ and hence also $T$ has range of dimension at most 
one, which is not the case. We infer that  $\deg T(h)=m$ for any 
$\overline{C}$-stable polynomial $h$ of degree $n$.

Clearly, $\phi_{m,z}(f)(z,w)=T[\phi_{n,z}\big((z+w)^n \big)]$. Let 
$w_0 \in H\setminus\{-a/c\}$.  The only zero of the degree $n$ polynomial  
$$p_{w_0}(z):=\phi_{n,z}\big( (z+w_0)^n \big)=\big( (a+w_0c)z+b+w_0d\big)^n$$ 
is $\Phi^{-1}(-w_0)\in \overline{C}'$, so $p_{w_0}(z)$ 
 is $\overline{C}$-stable. By the previous paragraph we then have
$$\deg(\phi_{m,z}(f)(z,w_0))= \deg (T[p_{w_0}(z)])=m,$$ 
which verifies condition (b1) of Lemma~\ref{fifi}. 

Note next that coefficient of $w^n$ in the polynomial defined in (b) 
of Theorem~\ref{circularregion} is $T[(2acz+bc+ad)^n]$. 
We claim that the polynomial $(2acz+bc+ad)^n$ is $\overline{C}$-stable and 
of degree $n$. Since the image of any $\overline{C}$-stable polynomial of 
degree $n$ is of degree $m$ this would verify condition (b2) of 
Lemma~\ref{fifi}. To prove the claim note first that since $C$ is the 
exterior of a disk we have that $ac\neq 0$ so $\deg(2acz+bc+ad)^n=n$. Let 
$\zeta=-b/a$ and $\eta= -d/c$. Since $\Phi(\zeta)=0$ and 
$\Phi(\eta)=\infty$ we have that $\zeta, \eta \in \partial C$, which 
implies -- again by the assumption that $C$ is the exterior of a 
disk --  that the zero of $(2acz+bc+ad)^n$
$$
-\frac {bc+ad}{2ac} = \frac {\zeta + \eta}{2}, 
$$
is in $\overline{C}'$. Thus $(2acz+bc+ad)^n$ is $\overline{C}$-stable and 
of degree $n$, as required.
\end{proof}

\begin{proof}[Proof of Theorem~\ref{circle}]
Let $m = \max\{\deg T(f): f \in \CC_n[z]\}$. By Lemma~\ref{save} we have 
that  $T : \pi_n(\partial C) \rightarrow \pi(\partial C)\cup\{0\}$ if and 
only if 
$\phi_m^{-1}T\phi_n : \pi_n(\RR) \rightarrow \pi(\RR)\cup\{0\}$. Using this 
and Theorem~\ref{finitehypC} it is not difficult to verify the theorem in the 
case $\dim_{\bC} T(\CC_n[z]) \leq 2$. Let 
$f(z,w)=\phi_m^{-1}T\phi_n [(z+ w)^n]$. To settle the remaining cases note 
first that by Theorem~\ref{finitehypC} (c) and (d) we have that 
$T : \pi_n(\partial C) \rightarrow \pi(\partial C)\cup\{0\}$ if and only if  
$f(z,w)$ or $f(z,-w)$ is a complex multiple of a real stable polynomial.  
Now the condition  that $f(z,w)$ is a complex multiple of a real stable 
polynomial is equivalent to saying that $f(z,w)$ is both 
$H$-stable and $H^r$-stable, see \cite[Proposition 3]{BBS} 
(note that $H^r=-H$). By Theorem~\ref{circularregion} we know that 
$f(z,w)$ is $H$-stable if and only if the polynomial in (c) of 
Theorem~\ref{circle}, that is, $\phi_{n,w} \phi_{m,z}(f)(z,w)$, is 
$C$-stable. On the other hand, since 
$-cz+a$ is $H^r$-stable and $cz+d$ is $C^r$-stable (since $C^r$ is not the 
exterior of a disk) we also have that $f(z,w)$ is $H^r$-stable if and only if 
$\phi_{n,w} \phi_{m,z}(f)$ is $C^r$-stable. The proof of 
the fact that condition (d) in Theorem~\ref{circle} is equivalent to saying 
that $f(z,-w)$ is a complex multiple of a real stable polynomial follows in 
similar fashion. This completes the proof of the theorem.    
\end{proof}

Corollary~\ref{alg-circularregion} and Corollary~\ref{cor-circle} are 
immediate consequences of Theorem~\ref{circularregion} and 
Theorem~\ref{circle}, respectively. 

\section{Open Problems}\label{s4}

As we already noted in \S \ref{s1}, Problems \ref{pb1} and \ref{pb2} have a
long and distinguished history. In this paper we completely solved them for
a particularly relevant type of sets, namely all closed circular domains and 
their boundaries. Among the most interesting remaining cases that are 
currently under investigation we mention:
\begin{enumerate}
\item[(a)] $\Omega$ is an open circular domain, 
\item[(b)] $\Omega$ is a sector or a double sector, 
\item[(c)] $\Omega$ is a strip, 
\item[(d)] $\Omega$ is a half-line,
\item[(e)] $\Om$ is an interval. 
\end{enumerate}

Let us briefly comment on the importance of the above cases. 

Case (a). The classical notion of  
Hurwitz (or continuous-time) stability refers to univariate polynomials with
all their roots in the open left half-plane. Its well-known discrete-time 
version -- often called Schur or Schur-Cohn stability -- is when all the 
roots of a 
polynomial lie in the open unit disk. Both these notions are fundamental and 
widely used in e.g.~control theory and engineering sciences. (The authors 
easily found  several hundreds of publications in both purely mathematical 
and applied areas devoted to the study of Hurwitz and Schur stability for 
various classes of polynomials as well as continuous or discrete-time 
systems.)  

Case (b). Polynomials and transcendental entire functions with all their zeros 
confined to a (double) sector often appear as solutions to 
e.g.~Schr\"odinger-type equations with polynomial potential or, more 
generally, any linear ordinary differential equation with polynomial 
coefficients and constant leading coefficient, see for instance 
\cite{HeRo, Hil}. 
Concrete information about linear transformations preserving this class of 
polynomials and entire functions turns out to be very useful for asymptotic 
integration of linear differential equations. 

Case (c). Specific examples of linear transformations preserving the 
class of polynomials with all their roots in the strip $|\Im(z)|\le \al$, 
$\al>0$, can be found in the famous articles by P\'olya \cite{Pol}, 
Lee-Yang \cite{LY} and de Bruijn \cite{deBr}. A complete characterization of 
all such linear transformations would shed light on a great many problems in 
the theory of Fourier transforms and operator theory \cite{BBCV,BBS2}.

Case (d). Polynomials with all their roots on either the positive or negative 
half-axis appeared already in Laguerre's works and in P\'olya-Schur's 
fundamental paper \cite{PS} and 
have been frequently used in various contexts 
ever since. In addition to their description of multiplier sequences of the
first kind in {\em op.~cit.~}the 
authors also classified multiplier sequences of the second kind, i.e., 
diagonal linear operators mapping polynomials with all real roots and of the 
same sign to polynomials with all real roots. A natural extension of these 
results would be to characterize all linear transformations with this
property. Solving case (d) would 
answer this question and thus complete the program initiated by P\'olya and 
Schur over 90 years ago. 

Case (e). Numerous papers have been devoted to this case of Problems~\ref{pb1}
and~\ref{pb2} and its connections with P\'olya frequency functions, integral
equations with totally positive or sign regular kernels, Laplace transforms, 
the theory of orthogonal polynomials, etc. A complete description of all 
linear transformations
preserving the set of polynomials with all their zeros in a given interval
would therefore have many interesting applications and would also answer 
several of the questions raised in 
e.g.~\cite{CPP2,iserles1,iserles2,iserles3,iserles4,Pin,PSz} on these 
and related subjects.

\begin{ack}
It is a pleasure to thank the American Institute of Mathematics for 
hosting the ``P\'olya-Schur-Lax Workshop'' \cite{BBCV} on topics pertaining
to this project in Spring 2007. 
\end{ack}

\end{document}